# Smarandache's function applied to perfect numbers

Sebastián Martín Ruiz
Avda. de Regla 43, Chipiona 11550 Spain
smruiz@telefonica.net
5 August 1998

**Abstract:**
Smarandache's function may be defined as follows:
S(n)= the smallest positive integer such that S(n)! is divisible by n. [1]

In this article we are going to see that the value this function takes when n is a perfect number of the form $n = 2^{k-1} \cdot (2^k - 1)$, $p = 2^k - 1$ being a prime number.

**Lemma 1**: Let $n = 2^i \cdot p$ when $P$ is an odd prime number and $i$ an integer such that:

$$0 \leq i \leq E\left(\frac{p}{2}\right) + E\left(\frac{p}{2^2}\right) + E\left(\frac{p}{2^3}\right) + \cdots + E\left(\frac{p}{2^{E(\log_2 p)}}\right) = e_2(p!)$$

where $e_2(p!)$ is the exponent of 2 in the prime number decomposition of $p!$.
E(x) is the greatest integer less than or equal to x.

One has that $S(n) = p$.

Demonstration:
Given that $GCD(2^i, p) = 1$ (GCD= greatest common divisor) one has that
$S(n) = \max\{S(2^i), S(p)\} \geq S(p) = p$. Therefore $S(n) \geq p$.
If we prove that p! is divisible by n then one would have the equality.

$$p! = p_1^{e_{p_1}(p!)} \cdot p_2^{e_{p_2}(p!)} \cdots p_s^{e_{p_s}(p!)}$$

where $p_i$ is the $i-th$ prime of the prime number decomposition of $p!$. It is clear that $p_1 = 2$, $p_s = p$, $e_{p_s}(p!) = 1$ for which:

$$p! = 2^{e_2(p!)} \cdot p_2^{e_{p_2}(p!)} \cdots p_{s-1}^{e_{p_{s-1}}(p!)} \cdot p$$

From where one can deduce that:

$$\frac{p!}{n} = 2^{e_2(p!)-i} \cdot p_2^{e_{p_2}(p!)} \cdots p_{s-1}^{e_{p_{s-1}}(p!)}$$

is a positive integer since $e_2(p!) - i \geq 0$.

Therefore one has that $S(n) = p$

**Proposition:** If n is a perfect number of the form $n = 2^{k-1} \cdot (2^k - 1)$ with $k$ is a positive integer, $2^k - 1 = p$ prime, one has that $S(n) = p$.

**Proof:**

For the Lemma it is sufficient to prove that $k-1 \leq e_2(p!)$.
If we can prove that:

$$k - 1 \leq 2^{k-1} - \frac{1}{2} \quad (1)$$

we will have proof of the proposition since:

$$k - 1 \leq 2^{k-1} - \frac{1}{2} = \frac{2^k - 1}{2} = \frac{p}{2}$$

As $k-1$ is an integer one has that $k - 1 \leq E\left(\frac{p}{2}\right) \leq e_2(p!)$

Proving (1) is the same as proving $k \leq 2^{k-1} + \frac{1}{2}$ at the same time, since k is integer, is equivalent to proving $k \leq 2^{k-1}$ (2).

In order to prove (2) we may consider the function: $f(x) = 2^{x-1} - x$ $\quad x$ real number.

This function may be derived and its derivate is $f'(x) = 2^{x-1} \ln 2 - 1$.

$f$ will be increasing when $2^{x-1} \ln 2 - 1 > 0$ resolving x:

$$x > 1 - \frac{\ln(\ln 2)}{\ln 2} \cong 1'5287$$

In particular $f$ will be increasing $\forall\ x \geq 2$.

Therefore $\forall\ x \geq 2 \quad f(x) \geq f(2) = 0$ that is to say $2^{x-1} - x \geq 0 \quad \forall x \geq 2$.

Therefore: $2^{k-1} \geq k \ \forall k \geq 2$ integer.

And thus is proved the proposition.

**EXAMPLES:**

$$\begin{array}{ll}
6 = 2 \cdot 3 & S(6)=3 \\
28 = 2^2 \cdot 7 & S(28)=7 \\
496 = 2^4 \cdot 31 & S(496)=31 \\
8128 = 2^6 \cdot 127 & S(8128)=127
\end{array}$$